\documentclass[12pt]{amsart}
\usepackage{amssymb}
\usepackage{amsxtra}
\usepackage[mathscr]{eucal}
\usepackage{graphics}
\usepackage{epsfig}

\renewcommand{\O}{\mathcal{O}}
\newtheorem{theorem}{Theorem}[section]
\newtheorem{lemma}[theorem]{Lemma}
\newtheorem{corollary}[theorem]{Corollary}
\newtheorem{proposition}[theorem]{Proposition}

\theoremstyle{definition}

\newtheorem{definition}[theorem]{Definition}

\newtheorem{example}[theorem]{Example}

\begin{document}
\title{Divisibility in elliptic divisibility sequences}
\author{Graham Everest and Helen King}
\address{(GE)School of Mathematics, University of East Anglia,
Norwich NR4 7TJ, UK.}\email{g.everest@uea.ac.uk}\email{h.king@uea.ac.uk}

\noindent
\dedicatory{\today}
\begin{abstract}
Certain elliptic divisibility sequences are shown to contain only
finitely many prime power terms. In certain circumstances, the
methods prove that only finitely many terms have length below a
given bound.
\end{abstract}
\maketitle

\section{Introduction}\label{intro}

Let $E$ denote an elliptic curve which is defined over $\mathbb
Q$. See \cite{cassels} or \cite{MR87g:11070} for background on
elliptic curves. The curve $E$ is given by a Weierstrass equation
\begin{equation}\label{weierform}
E:\quad y^2+a_1xy+a_3y=x^3+a_2x^2+a_4x+a_6,
\end{equation}
with $a_1,\dots a_6 \in \mathbb Z$. Suppose $E$ has a non-torsion
point $Q\in E(\mathbb Q)$. Throughout the paper, when $Q$ denotes
a point on an elliptic curve given in the form (\ref{weierform}),
the expression $x(Q)$ will denote the $x$-coordinate. For any
non-zero $n\in \mathbb Z$, write
\begin{equation}\label{AnBn}
x(nQ)=\frac{A_n}{B_n^2},
\end{equation}
in lowest terms, with $A_n$ and $B_n$ in $\mathbb Z$. The sequence
$B$ is a divisibility sequence, meaning that $B_m|B_n$ whenever
$m|n$ (a proof is supplied later). Since the sequence arises from
an elliptic curve, it seems natural to call $B$ an {\it elliptic
divisibility sequence}. In the literature (see \cite{MR9:332j} for
example) this term has been used for a related sequence, which
will be discussed at the end of the introduction. We agree with
Silverman \cite{silgcd} that the term is better used in the sense
proposed now.

In \cite{chuds}, Chudnovsky and Chudnovsky considered the
likelihood that the sequences $B$ might be a source of `large'
primes. They gave some impressive examples of prime values of
$B_n$.
\begin{example}
$$E:\quad y^2=x^3+26,\quad Q=[-1,5]
$$
The term $B_{29}$ is a prime with 285 decimal digits.
$$E:\quad y^2=x^3+15,\quad Q=[1,4]
$$
The term $B_{41}$ is a prime with 509 decimal digits.
\end{example}
In \cite{eds}, further numerical work suggested that for fixed $E$
and $Q$, the sequence $B_n$ should only contain finitely many
primes. A probabilistic argument, coupled with an affirmative
answer to the elliptic Lehmer problem, suggested the number of
prime terms should be uniformly bounded, see \cite{eds},
\cite{recseq}.

An {\it isogeny} between elliptic curves is a group homomorphism
defined by rational functions on the coordinates. Such a map has a
degree $d$, essentially the degree of the polynomials involved in
the rational map, and we say the map is a $d$-isogeny. The most
basic example is multiplication $Q\mapsto mQ$ where $m$ is an
integer. The degree of this isogeny is $m^2$. Any such map
factorizes as a composition of two isogenies, each of degree $m$.
In \cite{primeds}, it was proved that only finitely many prime
power terms $B_n$ appear if $Q$ is the image under a non-trivial
isogeny of a rational point. In \cite{primeds} we defined such a
point $Q$ to be {\it magnified}\footnote{The term was chosen
because the height of a point increases under such a map.}. This
definition will now be generalized.

Suppose $K$ denotes an algebraic number field and $Q\in E(K)$ is
non-torsion. Let $L/K$ denote an algebraic extension of finite degree. Suppose for
some $d>1$ there exists a $d$-isogeny $\sigma:F\rightarrow E$, from
an elliptic curve $F$,
mapping $R\in F(L)$ onto $Q$. We say $Q$ is {\it magnified from
R} if $[L:K]<d$. This definition
clearly generalizes the one in \cite{primeds} because the
definition in that paper applied only when $L=K=\mathbb Q$.
Let the set of valuations on $K$ be denoted by $M_K$, similarly
for $L$. (See section \ref{hival} for the definition of valuation.)
In the appendix to the paper an alternative way
of generalizing the definition of magnified will be discussed.

\begin{definition}\label{defofS}For any number field $K$, and any
point $Q\in E(K)$ let $S_K(Q)$ denote
the non-archimedean valuations $v$ in $M_K$ for which $|x(Q)|_v>1$.
\end{definition}

\begin{theorem}\label{finite}
Let $E$ denote an elliptic curve and let $Q\in E(K)$ denote a
non-torsion $K$-rational point, magnified from a point generating
a Galois extension. Then there are only finitely many $n\in
\mathbb N$ for which $S_K(nQ)$ consists of a single valuation.
\end{theorem}

\begin{corollary}\label{cor}
Let $E$ denote an elliptic curve and let $Q\in E(\mathbb Q)$
denote a non-torsion rational point which is magnified from a
point generating a Galois extension. Let $B$ denote the associated
elliptic divisibility sequence as in (\ref{AnBn}). Then only
finitely many terms $B_n$ are prime powers.
\end{corollary}

Theorem \ref{finite} clearly generalizes the main result in
\cite{primeds} and it will be proved using different methods. We
expect Theorem \ref{finite} is true without the assumption about
$L/K$ being Galois. Without this assumption, it is still possible
to prove that only finitely many terms $B_n$ are primes if $Q$ is
a magnified point. Later we present many examples of rational
points magnified from points generating quadratic extensions so
the Galois condition comes free for such examples. Later we will
also prove that sometimes the presence of a point $R$ with $mR=Q$
generating a Galois extension can force $Q$ to be the image of a
rational point under an $m$-isogeny. Thus, in some circumstances
we recover the results in \cite{primeds}.

Our generalization of the result in \cite{primeds} will be applied
in two ways. The first application gives a hint to what might be
ultimately true. If $l$ denotes a positive integer we say the
point $Q\in E(K)$ is {\it l-magnified} if it is the image of a
point under $l$ magnifications. If an integer has $r$ distinct
prime factors we say it has {\it length} $r$.

\begin{theorem}\label{length}
Let $E$ denote an elliptic curve and let $Q\in E(\mathbb Q)$
denote a non-torsion, rational $l$-magnified point. Let $B$ denote
the associated elliptic divisibility sequence as in (\ref{AnBn}).
Then only finitely many terms $B_n$ have length bounded by $l$.
\end{theorem}

{\bf Example}

1. Consider the elliptic curve
$$y^2=x^3-x^2-4x-2.
$$
The point $Q=[3,2]$ lies on $E$. Let $a^2-4a-4=0$ then
either point $R$ with $x(R)=a$ satisifes $2R=Q$. Thus $Q$ is a
magnified
point. Now let $b^4-16b^3-24b^2-16b-8=0$ then either point
$S$ with $x(S)=b$
satisfies $2S=R$ hence $R$ is itself
magnified. In both cases the Galois condition is satisfied because
the extensions are quadratic. Theorem \ref{length} shows the
equation $B_n=p^eq^f$
with $p$ and $q$ both primes has only finitely many solutions.

{\bf Conjecture} Suppose $E$ denotes an elliptic curve and $Q$
denotes a non-torsion rational point on $E$. We conjecture that if
the length of $B_n$ is bounded then $n$ is bounded. Further, we
believe that provided the curve $E$ is in minimal form, the bound
depends upon the length only and not $Q$ or $E$. In the terms of
Theorem \ref{finite}, we expect that if $|S_K(nQ)|$ is bounded then
$n$ is bounded without any assumptions about $L/K$.

The probabilistic arguments used in \cite{eds} (see also
\cite{recseq}) apply to support this conjecture. On the other hand
it is not easy to obtain convincing numerical evidence. Efficient
tests can be applied for primality (at least for `probable
primality' in the usual sense of computational number theory).
However no such tests can be applied, for example, to decide if an
integer is divisible by only two primes. This appears to require
the factorization of the integer and, given the rapid growth rate
of the sequences $B$, is impractical.

The second application of Theorem \ref{finite} is the following.
For any elliptic curve, and any integer $k\in \mathbb Z$, there is
an algebraic point $Q$ with $x(Q)=k$. The $x$-coordinates $x(nQ)$
for $n\in \mathbb N$ are all rational (assuming $Q$ is
non-torsion) so an integral sequence $B$ can be defined, just as
in (\ref{AnBn}). This provides a much richer stock of examples;
even if the elliptic curve has rank 0 it can still give rise to
integral divisibility sequences $B$. Theorem \ref{finite} applies
in this situation and later examples are provided of (singly and
multiply) magnified points $Q$.

{\bf Example} Suppose $E$ denotes the elliptic curve
$$
y^2+xy+y=x^3-36x-70.
$$

The point $Q$ with $x(Q)=-5$ gives rise to an integral sequence
$B$ as in $(\ref{AnBn})$. It is doubly magnified in one step by
times 2 and in the other by times 3, from points generating
respectively a quadratic and a cubic extension. The Galois closure
of the second extension over the first has degree at most 6
($<9$). Theorem \ref{length} implies the equation $B_n=p^eq^f$
with $p$ and $q$ both primes bounds $n$.

{\bf Nomenclature}
Finally a word about the terms used in this paper. Morgan Ward
used the term {\it elliptic divisibility sequence} (see
\cite{shipsey-thesis}, \cite{swart-thesis} or the papers
\cite{eds}, \cite{MR9:332j}) for bi-sequences (mapping $\mathbb Z$
to $\mathbb Z$) which satisfy
$$u_{n+k}u_{n-k}=u_{n-1}u_{n+1}u_{k}^2-u_{k-1}u_{k+1}u_n^2
$$
for all $n$ and $k$. These sequences are intimately connected with
the theory of elliptic curves. Given $E$ and $Q$ we can associate
a division sequence $\psi_n(Q)$ using the division polynomials,
which satisfy the recurrence relation above, see \cite{eds} for
details. Assuming without loss of generality that $Q$ is integral,
$\psi_n(Q)$ is an integral sequence with $B_n|\psi_n(Q)$. It
follows that if only finitely many terms $B_n$ have length bounded
by $l$, the same result follows for $\psi_n(Q)$ a fortiori.

\section{Valuations and Heights}\label{hival}

The valuations on $\mathbb Q$ consist of the usual archimedean
absolute value together with the non-archimedean,
$p$-adic valuations, written $|.|_p$ -
one for
each prime $p$. Let $v$ denote a valuation on $K$ then either
$v$ corresponds to an embedding of $K$ into $\mathbb C$ or it
corresponds to a prime ideal $\wp$ which is a factor of a
rational prime. Again, these are called archimedean and
non-archimedean valuations. Write $K_v$
for the completion of $K$ with respect to $v$. If $v$ is archimedean,
corresponding to the embedding $\phi:K\rightarrow \mathbb C$ then
$$|x|_v=|\phi(x)|^{1/[K:\mathbb Q]}.
$$
If $v$ is non-archimedean then it corresponds to the prime ideal $\wp$
dividing the rational prime $p$. Now
$$|x|_v=|x|_{\wp}^{e_v/[K:\mathbb Q]},
$$
where $e_v=[K_v:\mathbb Q_p]$. Given any $0 \neq \alpha \in K$,
the product formula holds:
$$\prod_{v\in M_K}|\alpha|_v=1,
$$
where the product runs over all valuations for $K$, both archimedean
and non-archimedean.
Write
\begin{equation}\label{localheight}
h_v(\alpha)=\log \max \{ 1, |\alpha|_v\},
\end{equation}
for the local (logarithmic) height at $v$.
The na\" ive global logarithmic height of $Q$ is defined to be
$$h(\alpha)=\sum_{v\in M_K} h_v(\alpha)=\sum_{v\in M_K} \log \max \{ 1, |\alpha|_v\},
$$
the sum running over all the valuations of $K$. The global height
is insensitive to the field of definition of $\alpha$ in the sense
that we can replace $K$ by any extension field in which $\alpha$
lies and the result will be the same. If $Q$ denotes any
$K$-rational point on an elliptic curve, we write $h_v(Q)$ and
$h(Q)$ for the above, when $\alpha=x(Q)$. Usually in the
literature, the global height is further normalized by dividing by
$2$.

Suppose $Q$ denotes a $K$-rational point of $E$. The
theory of heights gives an estimate for
$$
h(Q)=\widehat h(Q)+O(1),
$$
where $\widehat h(Q)$ denotes the canonical height of $Q$. The
canonical height enjoys the additional property that $\hat h(mQ)=m^2\hat h(Q)$
for any $m\in \mathbb Z$. More generally, if $\sigma:F\rightarrow E$
denotes a $d$-isogeny then for $R\in F(L)$
\begin{equation}\label{goingup}\hat h(\sigma(R))=d\hat h(R),
\end{equation}

\begin{lemma}\label{growth}
Suppose $Q\in E(K)$ denotes a non-torsion $K$-rational point. Then
$$h_v(nQ)=O(\log n \log \log n),
$$
for any archimedean valuation $v$.
\end{lemma}
\begin{proof}[Proof of Lemma~\ref{growth}.]
The estimate in Lemma~\ref{growth} follows from an appropriate
upper bound for $|x(nQ)|_v$. The archimedean valuations correspond
to the embeddings of $K$ into $\mathbb C$. We assume such a
valuation has now been fixed and work with the usual absolute
value on $\mathbb C$. Putting the model (\ref{weierform}) into the
usual Weierstrass equation only translates $x$ by at most a
constant. Thus we may assume that the $x$-coordinate of a point is
given using the Weierstrass $\wp$-function with Laurent expansion
in even powers of $z$,
$$x=\wp_L(z)=\frac{1}{z^2}+c_0+c_2z^2+\dots
$$
Let $z_Q$ correspond to $Q$ under an analytic isomorphism
$E(\mathbb C)\simeq \mathbb C/L$, for some lattice $L$. Write
$\{nz_Q\}$ for $nz_Q$ modulo $L$. When the quantity $|x(nQ)|$ is
large it means $nz_Q$ is close to zero modulo $L$, thus the
quantities $|x(nQ)|$ and $1/|\{nz_Q\}|^2$ are commensurate. On the
complex torus, this means the elliptic logarithm is close to zero.
So it is sufficient to supply a lower bound for $\{ nz_Q\}$ and
this can be given by elliptic transcendence theory (see
\cite{dav}). We use Th\'eor\`eme 2.1 in~\cite{dav} but see
also~\cite{ST} where an explicit version of David's Theorem
appears on page 20. The nature of the bound is
\begin{equation}\label{b2}
\log |x(nQ)| \ll \log n \log \log n,
\end{equation}
where the implied constant depends upon $E$, the valuation $v$ and
the point $Q$.
\end{proof}

We will need some theory of elliptic curves over local
fields,
see~\cite{MR87g:11070}. For every non-archimedean valuation $w$,
write $L_{w}$ for the completion of $L$ with respect to $w$. Write
ord$_w$ for the corresponding order function defined in terms of
the prime ideal associated to $w$. There is a subgroup of the
group of $L_{w}$-rational points:
$$E_1(L_{w})=\{O\} \cup \{Q\in E(L_{w}):\mbox{ord} _{w}(x(Q))
\leq -2\}.$$ In \cite{MR87g:11070}, Silverman proves the
following.
\begin{proposition}For all $R\in F_1(L_{w})$ and all $m\in \mathbb Z$:
\begin{equation}\label{ordthing}
\log|x(mR)|_w=\log |x(R)|_w-\log |m|_w.
\end{equation}
\end{proposition}

This important result yields several corollaries. The first is a
simple consequence of (\ref{ordthing}).
\begin{corollary}\label{divseq2}
The sequence $B$ is a divisibility sequence, meaning that
$B_m|B_n$, whenever $m|n$.
\end{corollary}

For finitely many prime ideals $\wp$ the reduction of $E$ or $F$
is not an elliptic curve because the reduced curve is singular.
Write $S$ for the set of valuations in $M_L$ corresponding to all
such prime ideals.

\begin{corollary}\label{bigOthing} Suppose $\sigma:F\rightarrow E$
denotes an isogeny and $R\in F_1(L_w)$. If $w \notin S$ then
$h_w(\sigma(R))\ge h_w(R)$. The local heights are related by the
formula
\begin{equation}\label{newordthing}
h_w(\sigma(R))=h_w(R)+O(1),
\end{equation}
where the implied constant depends only upon the isogeny and is
independent of $R$.
\end{corollary}

\begin{proof}[Proof of Corollary \ref{bigOthing}]
Any isogeny factorizes as a composition of isogenies of prime
degree. Without loss of generality, assume $\sigma$ is an isogeny
of prime degree $m$. Suppose $w$ corresponds to the prime ideal
$\wp$. Provided $w \notin S$ both curves and the isogeny can be
reduced modulo powers of $\wp$ and the first statement in the
Lemma follows. Applying the dual isogeny gives a similar
inequality $h_w(\sigma^*(Q))\ge h_w(Q)$ for all $Q\in E_1(L_w)$.
However, composing $\sigma$ with its dual gives multiplication by
$m$ on $F$. Now (\ref{ordthing}) applies to prove
(\ref{newordthing}).
\end{proof}

\begin{corollary}\label{fixedset}
Suppose $S$ denotes any fixed, finite set of valuations of $L$.
Then
\begin{equation}\sum_{w\in S}h_w(nR)=O((\log n)^2).
\end{equation}
\end{corollary}

\begin{proof}This is an immediate consequence of Lemma
\ref{growth}, for the archimedean valuations in $S$ and
(\ref{ordthing}) for the non-archimedean valuations in $S$.
\end{proof}

\begin{proof}[Proof of Theorem \ref{finite}]
Suppose there is an $L$-rational point $R$ with $\sigma(R)=Q$
where $\sigma:F\rightarrow E$ denotes a $d$-isgoeny. From
(\ref{goingup}), $\hat h(Q)=d\hat h(R)$, where $\hat h$ denotes
the canonical height of the algebraic points on $E$ and $F$. The
canonical height differs from the naive height by a bounded amount
we are justified in using only the naive height in the sequel.
\begin{equation}\label{heightdiff}
O(1)=dh(nR)-h(nQ).
\end{equation}
Suppose $S_K(nQ)$ consists of the single valuation $v$. We may
assume $v\notin S$ by Corollary \ref{fixedset}. From Lemma
\ref{growth}
$$h(nQ)-h_v(nQ)=O((\log n)^2)=h(nR)-\sum_{w\in S_L(nR)}h_w(nR),
$$
recalling that $S_L(nR)$ denotes those $w|v$ with $h_w(nR)>0$. It
follows from (\ref{heightdiff}) that
$$O((\log n)^2)=d\sum_{w\in S_L(nR)}h_w(nR)-h_v(nQ).
$$

The right hand side is
\begin{equation}\label{thisdoesit1}
d\sum_{w\in S_L(nR)}h_w(nR)-\sum_{w\in S_L(nQ)}h_w(nQ).
\end{equation}
The Galois assumption implies that the local height $h_w(nR)$ is
the same for each $w\in S_L(nR)$ and similarly for $h_w(nQ)$. Fix
$w\in S_L(nR)\subset S_L(nQ)$ by the first part of Lemma
\ref{bigOthing}. Write $e=|S_L(nR)|$ and $f=|S_K(nQ)|$.
By the Galois
assumption,
\begin{equation}\label{dividingthing}
e \quad | \quad f \quad | \quad [L_w:K_v] \quad | \quad [L:K].
\end{equation}

Hence (\ref{thisdoesit1}) becomes
\begin{equation}\label{thisdoesit2}
deh_w(nR)-fh_w(nQ).
\end{equation}
Now (\ref{newordthing}) implies
$$h_w(nQ)=h_w(nR)+O(1).
$$
Thus (\ref{thisdoesit2}) can be written
\begin{equation}\label{thisdoesit3}
(de-f)h_w(nQ) +O(1).
\end{equation}
It follows from (\ref{dividingthing}) that $(de-f) > 0$ if
$d>[L:K]$. This is enough to bound $n$ because $h_w(nQ)$ is
quadratic in $n$ while the expression in (\ref{thisdoesit3}) is
meant to be $O((\log n)^2)$.
\end{proof}

\begin{proof}[Proof of Theorem \ref{length}]This follows by
induction
using the methods in the proof of Theorem \ref{finite}. Suppose
that $Q=\sigma(R)$ is magnified via a  $d$-isogeny $\sigma$,
where $R$ generates
a Galois extension $L$.
For a contradiction, assume the only valuations
in $S_K(nQ)$ consist of those which
are extended by those in $S_{L}(nR)$.
As before
the expression (\ref{thisdoesit1})
is $O((\log n)^2)$.
Suppose there are $g$ conjugacy classes of valuations
in $S_{L}(nR)$ and choose representatives $w_1,\dots ,w_g$.
The Galois assumption implies that the local height is
the same for each valuation in a fixed class.
Write $e_i$ and $f_i$, as before, corresponding to each fixed class.
Then, by the Galois
assumption,
\begin{equation}\label{dividingthing2}
e_i \quad | \quad f_i \quad | \quad [L:K].
\end{equation}
Hence (\ref{thisdoesit1}) becomes
\begin{equation}\label{thisdoesitgen2}
\sum_{i=1}^g[de_ih_{w_i}(nR)-f_ih_{w_i}(nQ)].
\end{equation}
As before, using (\ref{newordthing}) allows
(\ref{thisdoesitgen2}) to be written
\begin{equation}\label{thisdoesitgen3}
\sum_{i=1}^g(de_i-f_i)h_{w_i}(nQ) +O(1).
\end{equation}
From (\ref{dividingthing2}), $(de_i-f_i) > 0$ if $d>[L:K]$. As before,
this is
enough to bound $n$ because the
expression in (\ref{thisdoesitgen3}) is meant to be $O((\log n)^2)$.
This contradiction shows that each time a point is magnified
under an isogeny, a new class of non-archimedean valuations ultimately
appears in the computation of the height.
\end{proof}

\section{Examples of magnified points}

Recall Velu's formulae for isogenies \cite{velu}. Every isogeny
factorizes as a product of isogenies of prime degree so consider
$m$ to be prime in what follows. Fix notation in the following
way: an isogeny $\sigma:F\rightarrow E$ will be described, where
$F$ is given in Weierstrass form
$$F:y^2+a_1'xy+a_3'y=x^3+a_2'x^2+a_4'x+a_6'.
$$
Firstly when $m=2$ there must be a rational 2-torsion point
$T=[x_1,y_1]$ on $F$. Write
\begin{eqnarray*}
t &=&3x_{1}^{2}+2a_{2}'x_{1}+a_{4}'-a_{1}'y_{1}, \\
u &=&4x_{1}^{3}+b_{2}'x_{1}^{2}+2b_{4}'x_{1}+b_{6}'\text{\qquad }\qquad \\
w &=&u+x_{1}t.
\end{eqnarray*}
Of course $u=0$ - this is simply a device to unify the
presentation. Then there is a 2-isogeny $\sigma$ to a curve $E$
taking $R=[x,y]$ to $\sigma(R)=Q$ where
$$x(Q)=x+t/(x-x_1).
$$
The isogenous curve $E$ has the form
\begin{equation}\label{formulaforisogenouscurve}
\lbrack a_{1},a_{2},a_{3},a_{4},a_{6}]=
[a_{1}',a_{2,}'a_{3,}'a_{4}'-5t,a_{6}'-b_{2}'t-7w].
\end{equation}
For odd $m,$ let $T=[x_{1},y_{1}]$ denote a point of order $m$ on
the curve and write its multiples as $kT=(x_{k},y_{k}),$ for
$1<k<m.$ The point $T$ may or may not have coordinates in the base
field $\mathbb Q.$ However the formulae, when applied to a point
$R$ on a curve $F$ over $\mathbb Q$, yield an isogenous curve $E$
and a point $Q$. The formulae for $t,u,w$ are given as follows.
Define, for each $k$ with $1\leq k\leq (m-1)/2,$
\begin{eqnarray*}
t_{k} &=&6x_{k}^{2}+b_{2}'x_{k}+b_{4}',\text{\qquad } \\
u_{k} &=&4x_{k}^{3}+b_{2}'x_{k}^{2}+2b_{4}'x_{k}+b_{6}'\qquad
\text{and\qquad }
\\
w_{k} &=&u_{k}+x_{k}t_{k}.
\end{eqnarray*}
Then with
\begin{equation*}
t=\sum t_{k}\qquad u=\sum u_{k}\qquad \text{and\qquad }w=\sum w_{k}
\end{equation*}
the formula for $E$ is exactly as in
(\ref{formulaforisogenouscurve}) and $x(Q)$ is given by
\begin{equation}\label{x-coordunderisog}
x(Q)=x+\sum_{k=1}^{(m-1)/2}\left\{
\frac{t_{k}}{(x-x_{k})}+\frac{u_{k} }{(x-x_{k})^{2}}\right\}.
\end{equation}

\begin{proposition}\label{2rat}Suppose $Q\in E(K)$ is non-torsion
point and $mR=Q$
where $R$ generates a Galois extension $L/K$. Suppose that for
all the points $R'$, which are Galois conjugates of $R$,
$R-R'$ is a $K$-rational $m$-torsion
point. Then $Q$ is the image of
a $K$-rational point under an isogeny.
\end{proposition}

When $m=2$, the hypothesis about $R-R'$ is always satisfied. This
is because $R+R'$ is a $K$-rational point which doubles to $2Q$
hence it differs from $Q$ by a $K$-rational torsion point. Write
$R+R'+T=Q$ then $T$ is non-trivial because if $R+R'=Q=2R$ then
$R=R'$ and $R$ is $K$-rational. Finally $R+R'+T=Q=2R$ implies
$R-R'=T$.

\begin{proof}
If $R$ is $K$-rational we are done. Otherwise, let $R'$ denote any
conjugate of $R$ not equal to $R$. Then $R-R'=T$ is a $K$-rational
$m$-torsion point on the curve $F$. Use this point to construct an
$m$-isogeny to an elliptic curve defined over $\mathbb Q$ as
above. Plainly $\sigma(T)=\O$. Hence
$\O=\sigma(R-R')=\sigma(R)-\sigma(R')=\sigma(R)-\sigma(R)'$ which
shows that $\sigma(R)\in E(K)$. Now the dual isogeny
$\sigma^*:F\rightarrow E$ has $\sigma ^* \sigma$ equal to
multiplication by $m$ on $E$. Hence $\sigma(R)$ is a $K$-rational
point which maps to $Q$ under the isogeny $\sigma^*$.
\end{proof}

Proposition \ref{2rat} is provable under even weaker hypotheses.
If all the points $R-R'$ can be written as multiples of a single
$m$-torsion point then that $m$-torsion point can be used in
Velu's formula to construct an $m$-isogeny $\sigma$ with the same
property as above - that $\sigma(R)$ is a $K$-rational point with
image $Q$ under the dual isogeny. Note however that the isogenous
curve is not necessarily defined over $K$.

The interest here is that we can construct points which are
multiply magnified. Table 1 show elliptic curves with multiply
magnified generators $Q$ under a multiplication by $m$ map. The
curves are recorded in the form $[a_1,a_2,a_3,a_4,a_6]$ to agree
with (\ref{weierform}). The examples were found among the first
500 rank-1 curves on Cremona's tables \cite{jec} using the
wonderful PARI-GP computing package \cite{pari}. For each curve,
we computed a factorization of the polynomial of degree $m^2$
whose roots are the $x$-coordinates of the points $R$ with $mR=Q$.
If the order of magnification is 2, there is a point $S$ with
$mS=R$, which generates a field of degree $m^2$. Order 3 indicates
the point $S$ is itself magnified from a point which generates a
field of degree $m^3$. Note that when $m=3$, Theorem \ref{finite}
is satisfied because the Galois closure of each extension has
degree at most 6.

Table 2 shows algebraic points $Q$ with $x(Q)=k$ associated to the
first few elliptic curves by conductor. In every case we searched
for $k=x(Q)$ with $|k|\le 20$. The table shows the curve $E$
together with $x$-coordinates of algebraic points magnified from
another algebraic point by a multiplication by $m$ map. Line 3 is
repeated at line 13 because the point shown is doubly magnified
via multiplication by 2 and 3.

{\bf Appendix} All of the main conclusions in the paper can be
given under an alternative generalization of the definition of the
term {\it magnified}. We could define the term to mean that $Q$ is
the image of an $L$-rational point under an isogeny of degree $d$
where $d$ is coprime to $[L:K]$. This condition can certainly be
fulfilled. On the other hand, in all the examples we have found, a
magnified point according to this definition is simultaneously
magnified by a rational point. We believe the definition could
have some value - a magnification in this sense could form part of
a chain of maps comprising a multiple magnification.

Table 3 shows rank-1 elliptic curves $E$ together with generators
$Q$. In each there is a 3-isogeny which maps a point $R$ to $Q$
where $R$ generates a quadratic extension. Besides this, there is
a 3-isogeny that maps a rational point to $Q$. These were computed
using Velu's formulae as above. The table also shows the order of
the rational torsion on $E$.

\begin{center}{\bf Table 1}
\begin{tabular}{|c|c|c|c|c|}
\hline
$m$ & $E$ & $Q$ & Torsion Order & Magnification\\
\hline
2 & $[1, -1, 1, 4, 6]$ & [0, 2] & 4 & 2\\
& $[0, -1, 0, -4, -2]$ & [3, 2] & 2 & 2\\
& $[1, -1, 1, 20, 22]$ & [8, 21] & 4 & 3\\
& $[1, -1, 1, -23, -34]$ & [-2, 1] & 2 & 3\\
& $[1, -1, 1, -37, 124]$ & [2, -9] & 4 & 2\\
& $[1, -1, 1, -67, 226]$ & [14, 36] & 4 & 2\\
& $[1, 1, 0, 4, 0]$ & [1, 2] & 2 &2\\
& $[1, 0, 0, -10, -13]$ & [7, 13] & 2 &2\\
& $[1, -1, 1, 7, -8]$ & [2, 2] & 2 &2\\
& $[1, 1, 1, -63, 156]$ & [6, 3] & 4 &2\\
& $[[1, 0, 1, 12, -14]$ & [2, 3] & 2 &2\\
& $[0, 0, 0, 9, 0]$ & [4, 10] & 2 &2\\
& $[1, -1, 0, -990, -11745]$ & [238, 3513] & 2 &3\\
& $[0, 1, 0, -39, -108]$ & [12, 36] & 2 &2\\
& $[1, -1, 0, 90, 436]$ & [12, 50] & 2 &2\\
\hline
3 & $[0, 1, 1, -7, 5]$ & [-1, 3] & 3 & 2\\
 & $[0, -1, 1, -65, -204]$ & [12, 24] & $1$ & 2\\
\hline
\end{tabular}
\end{center}

\bigskip
\begin{center}{\bf Table 2}
\begin{tabular}{|c|c|c|c|c|}
\hline
$m$ & $E$ & $k$ & Torsion Order & Magnification\\
\hline
2&$[1,0,1,-36,-70]$ & 8 & 6 & 1\\
&$$& 7 &  & 1\\
&$$ & -5 &  & 2\\
\hline
&$[1,0,1,-171,-874]$ & -8 & 2 & 1\\
&$[1,0,1,-11,12]$ & -4 & 6 & 1\\
&$[1,1,1,-10,-10]$ & -4 & 8 & 2\\
\hline
3&$[1,0,1,4,-6]$ & -1 & 6 & 1 \\
& & -3 &  & 1 \\
& & -19 &  & 1 \\
\hline
3&$[1,0,1,-36,-70]$ & 9 & 6 & 1\\
&$$ & 4 &  & 1\\
&$$ & -5 &  & 2\\
&$$ & -1 &  & 1\\
&$$ & -3 &  & 1\\
&$$ & -9 &  & 1\\
\hline
3&$[1,0,1,-1,0]$ & -9 & 6 & 2\\
\hline
3&$[1,0,1,-11,12]$ & -5 & 6 & 1\\
&$$ & -7 &  & 1\\
\hline
4&$[1,1,1,-10,-10]$ & -13 & 8 & 1\\
\hline
5&$[0, -1, 1, 0,0]$ & -1 & 5 &1\\
\hline
\end{tabular}
\end{center}

\begin{center}
{\bf Table 3}
\end{center}
\begin{center}
\begin{tabular}{|c|c|c|}
\hline
$E$ & $Q$ & Torsion Order \\
\hline
$[0, 1, 1, -7, 5]$&[-1, 3] &3 \\
$[0, 0, 1, -24, 45]$&[-3, 9] &3 \\
$[0, -1, 0, 8, -16]$&[4, 8] &1 \\
$[0, 1, 1, -27, 55]$&[-5, 9] &3 \\
$[1, 0, 0, -2, 4]$ &[-2, 2] &3\\
$[0, 1, 1, 3, 2]$ &[2, 4] &3\\
$[1, 0, 0, 3, 1]$ &[0, 1] &1 \\
$[0, 1, 0, -5, -5]$ &[-2, 1] &2 \\
$[0, -1, 0, 2, -1]$ &[1, 1] &1 \\
$[1, 0, 0, 4, 16]$ &[-2, 0] &3 \\
$[1, 1, 0, 3, 1]$ &[0, 1] &1 \\
$[1, -1, 1, 46, 209]$ &[-3, 7] &3 \\
$[0, 0, 1, -42, 105]$ &[5, 4] &3 \\
$[0, -1, 0, -53, -131]$ &[-4, 1] &1\\
$[0, 1, 0, -101, 359]$ &[2, 13] &3 \\
\hline
\end{tabular}
\end{center}


\end{document}